\documentclass[12pt]{article}

\usepackage{amsthm,amsmath,amssymb,fancyhdr}

\usepackage{graphicx}

\usepackage[colorlinks=true,citecolor=black,linkcolor=black,urlcolor=blue]{hyperref}

\theoremstyle{plain}
\newtheorem{theorem}{Theorem}
\newtheorem{lemma}[theorem]{Lemma}
\newtheorem{corollary}[theorem]{Corollary}
\newtheorem{proposition}[theorem]{Proposition}

\theoremstyle{definition}
\newtheorem{definition}[theorem]{Definition}

\theoremstyle{remark}
\newtheorem{remark}[theorem]{Remark}

\def\|{\Big |}
\def\({\Big (}
\def\){\Big )}
\def\[{\Big[}
\def\]{\Big]}
\def\be{\begin{equation}}
\def\bel{\begin{equation}\label}
\def\ee{\end{equation}}
\def\bt{\begin{theorem}}
\def\bcd{\begin{condition}}
\def\ecd{\end{condition}}
\def\et{\end{theorem}}
\def\bc{\begin{corollary}}
\def\ec{\end{corollary}}
\def\bde{\begin{definition}}
\def\ede{\end{definition}}
\def\bl{\begin{lemma}}
\def\el{\end{lemma}}
\def\bp{\begin{proposition}}
\def\ep{\end{proposition}}
\def\br{\begin{remark}}
\def\er{\end{remark}}
\def\ba{\begin{array}}
\def\ea{\end{array}}
\def\ed{\end{document}}
\def\Z{\mathbb Z}

\def\1{{\bf 1}}

\title{\bf On a sumset problem for integers}

\author{Shan-Shan Du\\
\small The Fundamental Division\\[-0.8ex]
\small Jinling Institute of Technology\\[-0.8ex]
\small 211169, People's Republic of China\\
\small\tt ssdu.stand@gmail.com\\
\and
Hui-Qin Cao\thanks{Supported by the National Natural Science Foundation (grant 11201233) of China.}\\
\small Department of Applied Mathematics\\[-0.8ex]
\small Nanjing Audit University\\[-0.8ex]
\small 211815, People's Republic of China\\
\small\tt caohq@nau.edu.cn\\
\and
Zhi-Wei Sun\thanks{Supported by the National Natural Science Foundation (grant 11171140) of China.}\\
\small Department of Mathematics\\[-0.8ex]
\small Nanjing
University\\[-0.8ex]
\small 210093, People's Republic of China\\
\small\tt zwsun@nju.edu.cn
}
\date{}

\begin{document}
\maketitle \thispagestyle{fancy} \lhead{Electron. J. Combin.
21(2014), no.1, \#P1.13, 1-25.}
\renewcommand{\headrulewidth}{0pt}
%
\begin{abstract}
  Let $A$ be a finite set
of integers. We show that if $k$ is a prime power or a product of
two distinct primes then
$$|A+k\cdot A|\geq(k+1)|A|-\lceil k(k+2)/4\rceil$$ provided $|A|\geq
(k-1)^{2}k!$, where $A+k\cdot A=\{a+kb:\ a,b\in A\}$. We also
establish the inequality $|A+4\cdot A|\geq5|A|-6 $ for $|A|\geq5$.

\end{abstract}

\section{Introduction}

 For finite subsets $A_1,\ldots,A_k$ of $\Z$, their {\it sumset} is given by
 $$A_1+\cdots+A_k=\{a_1+\cdots+a_k:\ a_1\in A_1,\ldots,a_k\in
 A_k\},$$
 which is simply denoted by $kA$ if $A_1=\cdots=A_k=A$.
 It is known that
$$|A_1+\cdots+A_k|\ge |A_1|+\cdots+|A_k|-k+1,$$
and equality holds when $A_1,\ldots,A_k$ are arithmetic progressions
with the same common difference (see, e.g., Nathanson
\cite[p.11]{n}).

 Let $A$ be a finite set of integers. For
 $k\in\Z^+=\{1,2,3,\ldots\}$, we define
 $$k\cdot A:=\{ka:a\in A\}$$
 which is called a {\it dilate} of $A$.
 Let $k_1,k_2,\ldots,k_l\in\Z^+$. Recently lower bounds for $|k_{1}\cdot A+k_{2}\cdot A+\cdots+k_{l}\cdot A|$ were investigated by various authors
 \cite{bb linear equation,bb sums of dilates,nooors,mbn 2008}. In the case $(k_{1}, k_{2},\ldots,k_{l})=1$
 (where $(k_1,\ldots,k_l)$ refers to the greatest common divisor of $k_1,\ldots,k_l$),
 Bukh \cite{bb sums of dilates} obtained the following inequality:
$$|k_{1}\cdot A+k_{2}\cdot A+\cdots+k_{l}\cdot A|\geq (k_{1}+k_2+\cdots+k_{l})|A|-o(|A|).$$
For $l=2$ there are better quantitative results in this direction, see \cite{chs,csv,hr,n,mbn 2008}. It was
conjectured in \cite{csv} that for any $k\in\Z^+$ if $|A|$ is sufficiently large then
 $$|A+k\cdot A|\geq(k+1)|A|-\lceil k(k+2)/4\rceil.$$
 This was proved in \cite{chs} with $k$ prime.
 In this paper we confirm the
 conjecture for $k=p^\alpha$ as well as $k=p_{1}p_{2}$, where $p$, $p_{1}$, $p_{2}$ are prime
 numbers and  $\alpha\in \mathbb{Z}^+$. Motivated by the preprint form of our paper posted to arXiv, Ljuji\'c \cite{zl}
 obtained similar results for $|2\cdot A+k\cdot A|$ with $k$ a prime power or a product of two distinct primes.

 We remark that there are also some researches on sums of dilates in $\mathbb Z_p=\Z/p\Z$ with $p$ a prime, see Plagne \cite{P}
 and Pontiveros \cite{Po}.

 Now we state our main theorems.

\begin{theorem} \label{t1} Let $k=p^\alpha$ with $p$  a prime and  $\alpha\in\Z^+$. Let $A$ be a finite subset of  $\mathbb{Z}$
with $|A|\geq (k-1)^{2}k!$. Then
 \be |A+k\cdot A|\geq(k+1)|A|-\lceil
k(k+2)/4\rceil.\ee
\end{theorem}

\begin{theorem}
\label{liangsushuchenji} Let $p_{1}$ and $p_{2}$ be distinct primes and
$k=p_{1}p_{2}$. And let $A$ be a finite subset of $\mathbb{Z}$ with
$|A|\geq (k-1)^{2}k!$. Then
\begin{equation}\label{2}
|A+k\cdot A|\geq(k+1)|A|-\lceil k(k+2)/4\rceil.
\end{equation}
\end{theorem}

By Theorem \ref{t1}, if $k=4$ then (1) holds when
$|A|\geq 216$. In fact, we have the following refinement.

\begin{theorem}
\label{k=4} For any finite set $A\subseteq\mathbb{Z}$ with $|A|\geq5$,
we have
\begin{equation}\label{3}
|A+4\cdot A|\geq5|A|-6.
\end{equation}
\end{theorem}

We remark that the lower bound given in (1) is optimal when $|A|$ is large enough. Moreover, equality holds if
 $A$ has the form $k\cdot\{0,1,\ldots,n\}+\{0,1,\ldots,h\}$, where
$$h= \left\{
  \begin{array}{ll}
  k/2 \ \mbox{or}\  (k+2)/2 &\ \mbox{if}\ k  \   \mbox{is even},\\
  (k+1)/2 & \ \mbox{if}\ k \  \mbox{is odd}.\\
  \end{array}
  \right.
$$

Our proofs of Theorems 1-3 are based on the technical approach of \cite{chs}. Our key new idea is to employ Chowla's theorem
to handle the case when $k$ is a prime power, and use a lemma similar to Chowla's theorem to handle the case when $k$ is a product of
two distinct primes.

\section{Preliminaries}

Throughout this
paper we use the following notations. For a finite set
$A\subseteq\mathbb{Z}$ with $|A|>1$ and a positive integer $k$, we
define $$\hat{A}=\{ \bar{a}=a+k\Z: a\in A\}.$$ Let $h=|\hat{A}|$ and
let $A_{1},A_{2},\ldots, A_{h}$ be the distinct classes of $A$
modulo $k$.  Write $A_{i}=k\cdot X_{i}+r_{i}$ with $0\leq r_{i}<k$ for
$i=1,2,\ldots,h$. Clearly $|A_{i}|=|X_{i}|$ and
$$A=\bigcup_{i=1}^hA_{i}=\bigcup_{i=1}^h(k\cdot X_{i}+r_{i}).$$
Define
$$F=\{1\leq i\leq h:|\hat{X_{i}}|=k\},\qquad  E=\{1\leq i\leq h:0<|\hat{X_{i}}|<k\}$$
and
$$\triangle _{rs}=(A_{r}+k\cdot A )\backslash(A_{r}+k\cdot
A_{s})\ \  \ \text{for} \ \ r,s=1,2,\ldots,h.$$

Without loss of generality, we make the following assumptions:
\vskip5pt
(I) $\gcd(A)=\gcd(\{a:a\in A\})=1$.

If $d=\gcd(A)>1$, then replace $A$ by
$A^{'}=\{a/d:a\in A\}$. Obviously $|A^{'}|=|A|$ and
$$|A^{'}+k\cdot A^{'}|=|A+k\cdot A|.$$

\medskip
(II) $r_{1}=0$ and $|A_{1}|\geq|A_{2}|\geq\cdots\geq|A_{h}|$.

In fact, for $A^{'}=A-r_{1}$ we have
$|A^{'}|=|A|$ and
$|A^{'}+k\cdot A^{'}|=|A+k\cdot A|$.
\medskip

(III) $h=|\hat{A}|\geq2$.

When $h=1$ we have $A=A_{1}=k\cdot X_{1}+r_{1}$ and $|A+k\cdot A|=|X_{1}+k\cdot X_{1}|$.
So we may replace $A$ by $X_{1}$, and continue this process until $|\hat{X_{1}}|>1$.

\begin{lemma}[cf. \cite{csv}] \label{a+kb} For arbitrary nonempty sets
$B$ and $A=\bigcup_{i=1}^h(k\cdot X_{i}+r_{i})$, we
 have\\
(i) \ \ $|A+k\cdot B|=\sum_{i=1}^h|X_{i}+B|$.\\
(ii)\ \  $|A+k\cdot B|\geq|A|+h(|B|-1)$.\\
(iii) Furthermore, if equality holds in (ii), then either $|B|=1$
 or $|X_{i}|=1$ for all $i=1,\ldots, h$ or $B$ and all the sets $X_{i}$ with more
than
 one element are arithmetic progressions with the same difference.
\end{lemma}

\begin{lemma}[I. Chowla, see \cite{n}]
\label{chowla} For $n\geq2$, let $A$ and $B$ be nonempty subsets
of
 $\mathbb{Z}/n\Z$. If $0\in B$ and $(b,n)=1$ for all $b\in
B\backslash\{0\}$, then
$$|A+B|\geq \min\{n, |A|+|B|-1\}.$$
\end{lemma}

\begin{lemma}[cf. \cite{chs}]
\label{ii} For each subset $I\subseteq{\{1,2,\ldots,h\}}$, we have
$$\sum_{i\in I}|\triangle_{ii}|\geq |I|(|I|-1).$$
\end{lemma}

\section{Proof of Theorem \ref{t1}}

In this section, we fix $k=p^\alpha$ where $p$ is a prime and $\alpha\in\mathbb{Z}^{+}$. Let $A$ be a nonempty finite subset of integers.
Note that the set $\{1\leq i\leq h:p\nmid r_{i}\}$ is nonempty since $\gcd (A)=1$. Define $m=\min \{1\leq i\leq
 h:p\nmid r_{i}\}$.
\begin{lemma}
\label{wo} Suppose that $A$ is a nonempty finite subset of integers.\\
(i) \  In the case  $i\in E\backslash \{m\}$, we have
$|\vartriangle_{ii}|\geq |A_{m}|$.\\
(ii) If $|\hat{X_{m}}|+m-1\leq k$, then
\begin{equation}\label{mm dayu 1,...m}
|\vartriangle_{mm}|\geq|A_{1}|+\ldots+|A_{m-1}|.
\end{equation}

\ If $|\hat{X_{m}}|+m-1> k$, then we have
\begin{equation}\label{am+a dang k wei sushumici}
|A_{m}+A|\geq(k+1)|A_{m}|+m(|A_{1}|-|A_{m}|)-k.
\end{equation}

\end{lemma}

\begin{proof}
(i)  Suppose $i\in E\backslash\{m\}$. Noting that $p\nmid r_{m}$, we
have $(r_{m}-r_{i}, k)=1$ when  $p\mid r_i$. Applying Lemma
\ref{chowla}, we get
$$
|\hat{X_{i}}+\{0, r_{m}-r_{i}\}|\geq\min\{k, |\hat{X_{i}}|+2-1\}=|\hat{X_{i}}|+1
$$
since $i\in E$. It follows that
 $$|(\hat{X_{i}}+\hat{A_{m}})\backslash(\hat{X_{i}}+\hat{A_{i}})|\geq1,$$ and hence,
\begin{align*}
  |\triangle _{ii}|&=|(A_{i}+k\cdot A )\backslash(A_{i}+k\cdot
  A_{i})|\\
  &=|(X_{i}+ A )\backslash(X_{i}+ A_{i})|\\
  &\geq|(X_{i}+ A_{m} )\backslash(X_{i}+ A_{i})|\\
  &\geq|(\hat{X_{i}}+\hat{A_{m}})\backslash(\hat{X_{i}}+\hat{A_{i}})|\cdot|A_{m}|\\
  &\geq|A_{m}|  \quad\mbox{(since}\  |\hat{A_{i}}|=|\hat{A_{m}}|=1).
\end{align*}
In the case $p\nmid r_{i}$, using Lemma \ref{chowla} we
obtain
\begin{align*}
|\hat{X_{i}}+\{0, r_{i}\}|&\geq \min\{k, |\hat{X_{i}}|+2-1\}=|\hat{X_{i}}|+1,
\end{align*}
hence
$|(\hat{X_{i}}+\hat{A_{1}})\backslash\hat{(X_{i}}+\hat{A_{i}})|\geq1$
and
\begin{align*}
  |\triangle _{ii}|&=|(A_{i}+k\cdot A )\backslash(A_{i}+k\cdot
  A_{i})|\geq|(X_{i}+ A_{1} )\backslash(X_{i}+ A_{i})|\geq|A_{1}|\geq|A_{m}|.
\end{align*}

(ii) \ Recall that $p\mid r_{1},\ldots,\ p\mid r_{m-1}$ and $p\nmid
r_{m}$. Thus
$$(r_{1}-r_{m},k)=\cdots=(r_{m-1}-r_{m},k)=1.$$ It follows from
Lemma \ref{chowla} that
$$|\hat{X_{m}}+\{0, r_{1}-r_{m}, r_{2}-r_{m},\ldots,r_{t}-r_{m}\}|\
\geq  \min \{k,\  |\hat{X_{m}}|+(t+1)-1|\}.$$ So we have
$$|\hat{X_{m}}+(\hat{A_{1}}\cup\cdots\cup\hat{A_{t}}\cup\hat{A_{m}})|
\geq  \min \{k,\  |\hat{X_{m}}|+t\}\ \ \text{for}\ \ t=1,2,\ldots,m-1.$$

 If $|\hat{X_{m}}|+m-1\leq k$, then by induction on $t$ we deduce that
 $$
|(X_m+A_1\cup\cdots\cup A_t\cup A_m)\backslash(X_{m}+
A_{m})|\geq|A_1|+\cdots+|A_t|
$$
for $t=1,2,\ldots,m-1$. Consequently,
\begin{align*}
  |\triangle _{mm}|&=|(X_{m}+ A )\backslash(X_{m}+ A_{m})|\geq|A_{1}|+\ldots+|A_{m-1}|.
\end{align*}

If $|\hat{X_{m}}|+m-1>k$, then
$|\hat{X_{m}}+(\hat{A_{1}}\cup\ldots\cup\hat{A_{m}})|=k$. With the
help of Lemmas \ref{a+kb} and \ref{chowla}, we get
\begin{align*}
  |X_{m}+A|&\geq |X_{m}+A_{1}|+|(X_m+A_1\cup\cdots\cup A_m)\backslash(X_{m}+ A_{1})|\\
  &\geq |X_{m}+k\cdot X_{1}|+|(\hat{X_{m}}+\hat{A_1}\cup \cdots \cup \hat{A_m})\backslash(\hat{X_{m}}+\hat{A_{1})}|\cdot|A_{m}|\\
  &\geq|X_{m}|+|\hat{X_{m}}|(|A_{1}|-1)+(k-|\hat{X_{m}}|)|A_{m}|\\
  &\geq(k+1)|A_{m}|+|\hat{X_{m}}|(|A_{1}|-|A_{m}|)-k.
\end{align*}
The definition of $m$ implies that $m\leq
p^{\alpha-1}+1$ and hence
\begin{align*}
  |\hat{X_{m}}|&> k+1-m\geq p^{\alpha}+1-p^{\alpha-1}-1\geq p^{\alpha-1}\geq m-1.
\end{align*}
Thus,
$$|X_{m}+A|\geq(k+1)|A_{m}|+m(|A_{1}|-|A_{m}|)-k.$$
\end{proof}

\begin{lemma}
\label{small} Let $A$ be a nonempty finite subset of  $\mathbb{Z}$.
Then
$$|A+k\cdot A|\geq(k+1)|A|-k!.$$
\end{lemma}

\begin{proof} It suffices to prove by induction that
\begin{equation}
\label{guinazhengming equation}
|A+k\cdot A|\geq(t+1)|A|-(t-1)!k.
\end{equation}
holds for every $t=1,\ldots,k$.

Clearly (\ref{guinazhengming equation}) is true for $t=1$ since it is
known that
$$|A+k\cdot A|\geq 2|A|-1\geq 2|A|-k.$$

Now suppose that (\ref{guinazhengming equation}) holds for some $1\le
t<k$. We want to deduce (\ref{guinazhengming equation}) with $t$
replaced by $t+1$, i.e., the inequality
$$
|A+k\cdot A|\geq (t+2)|A|-t!k.
$$
If $h>t$, then applying Lemma \ref{a+kb} we immediately get
$$
|A+k\cdot A|\geq |A|+h(|A|-1)\geq|A|+(t+1)(|A|-1)\geq(t+2)|A|-t!k.
$$
Below we assume $h\leq t$. By Lemma \ref{a+kb},  for $i\in F$ we
have
\begin{align*}
 |A_{i}+k\cdot A|&\geq|X_{i}+k\cdot X_{1}|\\
  &\geq |X_{i}|+k|X_{1}|-k\\
  &\geq |A_{i}|+(t+1)|A_{1}|-k\\
  &\geq(t+1)|A_{i}|+|A_{1}|-k.
\end{align*}
By the induction hypothesis and Lemma \ref{wo}, for $i\in
E\setminus\{m\}$ we get
\begin{align*}
  |A_{i}+k\cdot A|&=|A_{i}+k\cdot
  A_{i}|+|\triangle_{ii}|\geq(t+1)|A_{i}|-(t-1)!k+|A_m|.
  \end{align*}
Therefore,
\begin{align*}
 \sum\limits_{i\not=m}|A_i+k \cdot A|
  &=\sum\limits_{i\in F\backslash\{m\}}|A_{i}+k\cdot A|
   +\sum\limits_{i\in E\backslash\{m\}}|A_{i}+k\cdot A|\\
  &\geq(t+1)\sum\limits_{i\not=m}|A_i|+|F\backslash\{m\}||A_{1}|+|E\backslash\{m\}||A_{m}| \\
   &\ \ \ -\(k|F\backslash\{m\}|+|E\backslash\{m\}|(t-1)!k\).
\end{align*}

We divide the following discussion into two cases.

{\it Case} 1. $|\hat{X_{m}}|+m-1\leq k$.

In this case, by (\ref{mm dayu 1,...m}) and the induction hypothesis, we have
\begin{align*}
|A_m+k \cdot A|&=|A_m+k \cdot A_m|+|\triangle_{mm}|\\
 &\geq(t+1)|A_m|-(t-1)!k+|A_1|+\cdots+|A_{m-1}|.
\end{align*}
It follows that
\begin{align*}
 |A+k \cdot A|
  &=\sum\limits_{i\not=m}|A_{i}+k\cdot A|
   +|A_{m}+k\cdot A|\\
   &\geq(t+1)|A|+|F\backslash\{m\}||A_{1}|+|E\backslash\{m\}||A_{m}|+|A_{1}|+\ldots+|A_{m-1}| \\
   &\ \ \ -\(k|F\backslash\{m\}|+|E\backslash\{m\}|(t-1)!k+(t-1)!k\).
\end{align*}
Clearly,
$$|F\backslash\{m\}||A_{1}|+|E\backslash\{m\}||A_{m}|+|A_{1}|+\ldots+|A_{m-1}|\geq|A|$$
and
$$k|F\backslash\{m\}|+|E\backslash\{m\}|(t-1)!k+(t-1)!k\leq (t-1)!kh\leq t!k.$$
Hence
$$|A+k\cdot A|\geq(t+2)|A|-t!k.$$

{\it Case} 2. $|\hat{X_{m}}|+m-1>k$.

We obtain from (\ref{am+a dang k wei sushumici}) that
\begin{align*}
  \ \ \ \ |A+k \cdot A|&=\sum\limits_{i\not=m}|A_{i}+k\cdot A|
   +|A_{m}+k\cdot A|\\
   &\geq\sum\limits_{i\not=m}|A_{i}+k\cdot A|
    +(k+1)|A_{m}|+m(|A_{1}|-|A_{m}|)-k\\
   &\geq(t+1)|A|+|F\backslash\{m\}||A_{1}|+|E\backslash\{m\}||A_{m}|+|A_{m}|+m(|A_{1}|-|A_{m}|)\\
   &\ \ \ -\(k|F\backslash\{m\}|+|E\backslash\{m\}|(t-1)!k+k\).
\end{align*}
As $|A_{1}|\geq|A_{2}|\geq\cdots\geq|A_{h}|$, we have
\begin{align*}
&|F\backslash\{m\}||A_{1}|+|E\backslash\{m\}||A_{m}|+|A_{m}|+m(|A_{1}|-|A_{m}|)\\
\geq&(|F|+|E|)|A_{m}|+m(|A_{1}|-|A_{m}|)\\
=&h|A_{m}|+m(|A_{1}|-|A_{m}|)=m|A_{1}|+(h-m)|A_{m}|\\
\geq&|A_{1}|+|A_{2}|+\cdots+|A_{m}|+|A_{m+1}|+\cdots+|A_{h}|=|A|
\end{align*}
and
$$k|F\backslash\{m\}|+|E\backslash\{m\}|(t-1)!k+k\leq (t-1)!kh\leq t!k.$$
Consequently,
$$|A+k\cdot A|\geq(t+2)|A|-t!k$$
as desired. This concludes the induction step.
\end{proof}

\begin{proof}[Proof of Theorem~\ref{t1}]
 Now suppose $|A|\geq (k-1)^{2}k!$. When $h=k$, Lemma \ref{a+kb} shows $|A+k\cdot A|\geq(k+1)|A|-k$, which means that (1) is valid.
 Below we assume $h<k$, and thus $|A|\geq(k-1)^{2}k!\geq h^{2}k!$, from which we have $|A_{1}|\geq |A|/h\geq hk!.$

{\it Case} 1. $i\in F$ for all $1\leq i\leq h$.

Due to Lemmas \ref{a+kb} and \ref{ii} we conclude that
\begin{align*}
|A+k\cdot A|&=\sum\limits_{i=1}^{h}|A_{i}+k\cdot
A|=\sum\limits_{i=1}^{h}\(|X_{i}+k\cdot
X_{i}|+|\triangle_{ii}|\)\\
&\geq\sum\limits_{i=1}^{h}\(|X_{i}|+k(|X_{i}|-1)+|\triangle_{ii}|\)\\
&\geq(k+1)|A|-hk+h(h-1)\\
&=(k+1)|A|-h(k+1-h).
\end{align*}
If $k$ is odd then
$$h(k+1-h)\leq\frac{k+1}{2}\(k+1-\frac{k+1}{2}\)=\frac{(k+1)^{2}}{4}=\Big\lceil\frac{k(k+2)}{4}\Big\rceil;$$
if $k$ is even then
$$h(k+1-h)\leq\frac{k}{2}\(k+1-\frac{k}{2}\)=\frac{k(k+2)}{4}=\Big\lceil\frac{k(k+2)}{4}\Big\rceil.$$
Therefore,
$$|A+k\cdot A|\geq(k+1)|A|-\lceil k(k+2)/4\rceil. $$
{\it Case} 2. $m\in E$.

We have
$|\hat{X_m}+\hat{A_1}\cup\hat{A_m}|\geq|\hat{X_m}|+1$ from Lemme \ref{chowla} and hence
$$|\triangle_{mm}|=|(A_m+k\cdot A)\setminus(A_m+k\cdot A_m)|\geq|(X_m+A_1)\setminus(X_m+A_m)|\geq|A_{1}|.$$
Then using Lemma \ref{small} we conclude that
\begin{align*}
|A+k\cdot A|&=|A_{m}+k\cdot A|+\sum\limits_{i\neq m}|A_{i}+k\cdot A|\\
&=|A_{m}+k\cdot A_{m}|+|\triangle_{mm}|+\sum\limits_{i\neq m}|A_{i}+k\cdot A|\\
&\geq(k+1)|A_{m}|-k!+|A_{1}| +\sum\limits_{i\neq m}\((k+1)|A_{i}|-k!\)\\
&\geq(k+1)|A|-hk!+|A_{1}|.
\end{align*}
By the fact $|A_1|\geq hk!$, we have
$$|A+k\cdot A|\geq(k+1)|A|.$$

{\it Case} 3. $m\in F$ and there exists $s\neq m$ such that $s\in E$.

In this case, Lemma \ref{wo} implies $ |\triangle_{ss}|\geq |A_{m}|$.
Then applying Lemmas \ref{a+kb} and \ref{small}, we see that
\begin{align*}
|A+k\cdot A|&=|A_{m}+k\cdot A|+\sum\limits_{i\neq m}|A_{i}+k\cdot A|\\
&\geq |A_{m}+k\cdot A_1|+\sum\limits_{i\neq m}|A_{i}+k\cdot A_i|+|\triangle_{ss}|\\
&\geq|A_{m}|+k|A_{1}|-k+\sum\limits_{i\neq m}\((k+1)|A_{i}|-k!\)+|A_{m}|\\
&\geq(k+1)|A|-hk!+|A_{1}|\\
&\geq (k+1)|A|.
\end{align*}

In view of the above discussions we have completed the proof of Theorem \ref{t1}.
\end{proof}

\section{Proof of Theorem \ref{liangsushuchenji}}

\begin{lemma}
\label{bing} Let $k$ be a positive integer and let $A$ be a nonempty
subset of \\$\mathbb{Z}/k\Z=\{\bar{0},
\bar{1},\ldots\overline{k-1}\}$. For $
\alpha\in\{1,2,\ldots,k-1\}$, we have $A+\bar{\alpha}=A$ if and only if
$$A=\bigcup_{\beta\in I} \left( (k,\alpha)\cdot\left\{\bar{0},\bar{1},\cdots,
\bar{\frac{k}{(k,\alpha)}}-\bar{1}\right\}+\bar{\beta}\right)$$ for some nonempty set
$I\subseteq \{0, 1,\ldots,(k, \alpha )-1\}$
.
\end{lemma}

\begin{proof}
In the case $(k, \alpha)=1$, it is easy to get that $A+\bar{\alpha}=A$
holds if and only if $A=\mathbb{Z}/k\Z$, which yields the Lemma.
Below we assume $(k,\alpha)>1$.

Let $A=A^{0}\cup A^{1}\cup \ldots\cup A^{(k,\alpha)-1}$ with
$$A^{i}=\{\gamma\in A: \gamma\equiv \bar{i}\  \ (\rm mod\ (k,\alpha))\}.$$
Note that $A+\bar{\alpha}=A$ implies $A^{i}+\bar{\alpha}=A^{i}$. Then by the fact $\(\frac{k}{(k,\alpha)},\frac{\alpha}{(k,\alpha)}\)=1$,
 we obtain \\ $$A^{i}=(k,\alpha)\cdot\left\{\bar{0},\bar{1},\cdots,
\bar{\frac{k}{(k,\alpha)}}-\bar{1}\right\}+\bar{i} \
 \  or \ \ \  A^{i}=\emptyset.$$ Thus
$$A=\bigcup_{\beta\in I} \left( (k,\alpha)\cdot\left\{\bar{0},\bar{1},\cdots,
\bar{\frac{k}{(k,\alpha)}}-\bar{1}\right\}+\bar{\beta}\right)$$ for some nonempty set
$I\subseteq \{0, 1,\ldots,(k, \alpha
)-1\}$.

The sufficiency is obvious, and the claim follows.
\end{proof}

\begin{lemma} [cf. \cite{n}]
\label{etransform} Let $A$ and $B$ be nonempty subsets of the
abelian group $G$, and let $g$ be any element of $G$. Let
$(A(g),B(g))$ be the e-transform of the pair (A,B), defined by
$A(g)=A\cup (B+g)$ and $B(g)=B\cap(A-g)$. Then $$A(g)+B(g)\subseteq
A+B$$ and$$A(g)\setminus A=g+(B\setminus B(g)).$$ If $A$ and $B$ are
finite sets, then
$$|A(g)|+|B(g)|=|A|+|B|.$$ If $g\in A$ and $0\in B$, then $g\in A(g)$ and $0\in
B(g)$.
\end{lemma}

The following lemma is a variation of Lemma \ref{chowla}.

\begin{lemma}
\label{gaijinchowla} Let $k>2$ be a composite integer. And let $A$,$B$ be
nonempty subsets of $\mathbb{Z}/k\Z$ with $A\neq \mathbb{Z}/k\Z$.
Assume
 $\bar{0}\in B$ and $\bar{0}\neq\bar{q}\in B$ with $(q,k)\neq 1$. If $(b,k)=1$ for all $\bar{b}\in B\backslash \{\bar{0},\bar{q}\}$,
 then
$$|A+\{\bar{0},\bar{q}\}|\geq |A|+1 \Rightarrow |A+B|\geq \min\{k,|A|+|B|-1\}.$$
\end{lemma}
\begin{proof}
Obviously it is true in the case $|A|+|B|>k$. Now we suppose
$|A|+|B|\leq k$. It is easy to deduce that the lemma holds
for $|A|=1$ or $|B|\leq2$. Next suppose $|A|\geq 2$ and $|B|\geq 3$.
If the claim fails, then there exist sets $A$, $B$ such that
$|A+B|<|A|+|B|-1$. Choose the pair $(A, B)$ such that $|B|$ is the
smallest. Since $|B|\geq 3$, we have $\bar{b^{*}}\in
B\setminus\{\bar{0},\bar{q}\}$. Then $(b^{*}, k)=1$. Due to
$A\not=\mathbb{Z}/k\Z$, there exists $\bar{g}\in A$
such that $\bar{g}+\bar{b^{*}}\notin A$ by Lemma \ref{bing}. Applying the $g$-transform to the pair
$(A, B)$ we have
$$|A(\bar{g})+B(\bar{g})|<|A(\bar{g})|+|B(\bar{g})|-1$$ and $$|B(\bar{g})|<|B|.$$
If $\bar{q}\in B(\bar{g})$, then it contradicts the minimality of $|B|$.
If $\bar{q}\notin  B(\bar{g})$, then we have $|A(\bar{g})+B(\bar{g})|\geq |A(\bar{g})|+|B(\bar{g})|-1$ from Lemmas \ref{chowla}
and \ref{etransform} , which
is also a contradiction.

This completes the proof.
\end{proof}

From now on we fix $k=p_{1}p_{2}$ in this section with $p_{1},
p_{2}$ distinct prime numbers. Suppose that $A$ is a nonempty finite
subset of $\mathbb{Z}$. In the case $(r_{2},k)> 1$, we may suppose
$(r_{2},k)=p_{1}$ without loss of generality. Then denote
$$n=\min\{\ 1\leq i\leq h: p_{1}\nmid r_{i}\}.$$
\begin{lemma}
\label{changsui} Let $A$ be a nonempty finite subset of $\mathbb{Z}$.

(i) If $(r_{2},k)=1$, then $|\triangle_{22}|\geq |A_1|$ for $2\in E$
and $|\triangle_{ii}|\geq |A_2|$ for $i\in E\backslash \{2\}$.

(ii)Suppose $(r_{2},k)= p_1$. Then
$$
|\triangle_{11}|\geq|A_{2}|\quad \text{or}\quad p_{2}|A_{n}| \qquad
\text{if}\ \ 1\in E,
$$
$$
|\triangle_{ii}|\geq|A_{1}|\quad \text{or}\quad p_{2}|A_{n}| \qquad
\text{for}\ \ i\in E\cap\{2,3,\ldots,n-1\}
$$
and
$$
|\triangle_{ii}|\geq|A_{n}| \qquad \text{for}\ \ i\in
E\cap\{n+1,\ldots,h\}.
$$
When $n\in E$, we have $|\triangle_{nn}|\geq|A_{2}|$. Moreover,
$$
  |X_{n}+A|\geq \begin{cases}|A_{n}|+p_{1}|A_{1}|-k
  &if\ |\hat{X_{n}}|\geq p_{1}>p_{2}, \\
  |A_{n}|+p_{2}|A_{1}|-k
  &if\ |\hat{X_{n}}|\geq p_{2}>p_{1}, \\
 |A_{n}|+|\hat{X_{n}}|\cdot|A_{1}|+|A_{2}|+\ldots+|A_{l}|-k
  &if\ p_{1}\leq|\hat{X_{n}}|< p_{2},
\end{cases}$$
where $l=\min\{n-1, p_{2}+1-|\hat{X_{n}}|\}$, and
 $$ |\triangle_{nn}|\geq|A_{1}|+|A_{2}|+\ldots+|A_{n-1}|\quad\mbox{if}\
|\hat{X_{n}}|< p_{1}.$$
\end{lemma}

\begin{proof}
(i) Note that $r_{1}=0$ and $(r_{2},k)=1$. Applying Lemma
\ref{chowla} we get
$$
|\hat{X_{i}}+\{0, r_{2}\}|\geq \min\{k, |\hat{X_{i}}|+2-1\}
$$
and hence
$$
\hat{X_i}+\hat{A_1}\not=\hat{X_i}+\hat{A_2}\qquad\text{for all}\qquad
i\in E.
$$
So we have that $|\triangle_{11}|\geq |A_{2}|$ for $1\in E$ and
that $|\triangle_{22}|\geq |A_{1}|$ for $2\in E$.\\
For $i\in E\backslash \{1,2\}$, if $(r_{i}, k)=1$ then by Lemma
\ref{chowla} we have
$$\hat{X_{i}}+\hat{A_{i}}\neq \hat{X_{i}}+\hat{A_{1}}.$$
Now suppose $(r_{i}, k)\neq 1$. In the case  $(r_{i}-r_{2}, k)=1$,
we have
$$\hat{X_{i}}+\hat{A_{i}}\neq \hat{X_{i}}+\hat{A_{2}}.$$
If $(r_{i}-r_{2}, k)\neq 1$, then $(r_{i}, k)\neq(r_{i}-r_{2}, k)$,
and hence by Lemma \ref{bing} we obtain
$$\hat{X_i}+\hat{A_i}\not=\hat{X_i}+\hat{A_1}\quad\text{or}\quad\hat{X_i}+\hat{A_i}\not=\hat{X_i}+\hat{A_2}.$$
Consequently,
 $$|\triangle_{ii}|=|(X_i+A)\setminus(X_i+A_i)|\geq |A_{2}|.$$

(ii) Suppose $1\in E$. If
$\hat{X_1}+\hat{A_1}\not=\hat{X_1}+\hat{A_2}$ then
$$|\triangle_{11}|\geq|(X_i+A_{2})\setminus(X_i+A_1)|\geq |A_{2}|.$$
In the case $\hat{X_1}=\hat{X_1}+\hat{A_1}=\hat{X_1}+\hat{A_2}$, by
Lemma \ref{bing} there is a proper subset $I$ of
$\{0,1,\ldots,p_{1}-1\}$ such that
$$\hat{X_{1}}=\bigcup_{\beta\in I} \(
p_{1}\cdot\{\bar{0},\bar{1},\cdots, \overline{p_{2}}-\bar{1}\}+\bar{\beta}\)$$ since $p_{1}=(r_{2}, k)$ and $1\in E$. Recall that $p_{1}\nmid
r_{n}$, and thus we
 have $$|(\hat{X_{1}}+\hat{A_{n}})\backslash
(\hat{X_{1}}+\hat{A_{1}})|\geq p_{2}$$ because of $I\neq \{0, 1,\ldots,p_{1}-1\}$, from which we get $$|\triangle_{11}|\geq
p_{2}|A_{n}|.$$
\ \ \ Similarly, for $i\in E\cap\{2,\ldots,n-1\}$, we have $$|\triangle_{ii}|\geq|A_{1}|\quad\text{or}\quad
p_{2}|A_{n}|.$$
\ \ \ If $i>n$
 and $i\in E$, then we also have
 $$|\triangle_{ii}|\geq|A_{1}|\quad\text{or}\quad
p_{2}|A_{n}|\geq |A_{n}|$$  when $(r_{i}, k)=1$ or $p_{1}$.
In the case $p_{2}| r_{i}$, we have $(r_{i}-r_{2}, k)=1$, and hence $\hat{X_i}+\hat{A_2}\not=\hat{X_i}+\hat{A_i}$
by Lemma \ref{bing}. So
$|\triangle_{ii}|\geq
 |A_{2}|\geq|A_{n}|$.

Below we discuss $|\triangle_{nn}|$ and  $|X_{n}+A|$ for $n\in E$. Since $p_{1}\nmid r_{n}$ and
$k=p_{1}p_{2}$, we have $(r_{n}, k)=1$ or
$(r_{n}-r_{2},k)=1$. Therefore
$$\hat{X_n}+\hat{A_n}\not=\hat{X_n}+\hat{A_1}\quad\text{or}\quad\hat{X_n}+\hat{A_n}\not=\hat{X_n}+\hat{A_2},$$
which states
$$|\triangle_{nn}|=|(X_{n}+A)\backslash(X_{n}+A_{n})|\geq
|A_{2}|.$$

Moreover with the help of Lemma \ref{a+kb}, we get
$$|X_{n}+A|\geq|X_{n}+A_1|\geq|X_{n}|+|\hat{X_{n}}|(|A_{1}|-1)\geq
|X_{n}|+|\hat{X_{n}}||A_{1}|-k,$$
and hence the claim holds for the case $|\hat{X_{n}}|\geq p_{1}>p_{2}$ or
 $|\hat{X_{n}}|\geq p_{2}>p_{1}$.

Now we turn to the last two cases.

{\it Case} 1. $p_{1}\leq |\hat{X_{n}}|< p_{2}$.

 Since
$|\hat{X_{n}}|< p_{2}$, we have $|\{x\ ({\rm mod}\  p_{2}): x\in
X_{n}\}|\leq |\hat{X_{n}}|<p_{2}$. Observing that $$(p_{2},
r_{2})=(p_{2}, r_{3})=\ldots=(p_{2}, r_{n-1})=1$$ and that
$$|\{r_{i}\ ({\rm mod}\  p_{2}): 2\leq i\leq n-1\}|=n-2,$$ in light of
Lemma \ref{chowla} we get
$$|\hat{X_{n}}+(\hat{A_{1}}\cup\cdots\cup\hat{A_{t}})|
\geq \min\{p_2, |\hat{X_n}|+t-1\} \quad\text{for}\quad t=1, 2,\ldots, n-1.$$ Hence
\begin{align*}
  |(X_{n}+A)\setminus(X_n+A_1)|&\geq|(X_n+A_1\cup\cdots\cup A_{n-1})\setminus(X_n+A_1)|\\
  &\geq |A_{2}|+\cdots+|A_{l}|,
\end{align*}
where $l=\min\{n-1, p_{2}+1-|\hat{X_{n}}|\}$. Consequently,
\begin{align*}
  |X_{n}+A|&\geq|X_{n}+A_{1}|+|A_{2}|+\ldots+|A_{l}|\\
  &\geq|A_{n}|+|\hat{X_{n}}||A_{1}|+|A_{2}|+\cdots+|A_{l}|-k.
\end{align*}

{\it Case} 2. $|\hat{X_{n}}|< p_{1}$.

By the definition of $n$, we have $n\leq p_{2}+1$ and hence
$$
|\hat{X_{n}}|+n-1<p_{1}+p_{2}\leq p_{1}p_{2}=k.
$$
Recall that $p_1\mid r_i$ for $1\leq i\leq n-1$ and that $p_{1}\nmid r_{n}$.
If there exists $1\leq s\leq n-1$ with $(r_{s}-r_{n}, k)\neq1$,  then we
 have $p_{2}\mid (r_{s}-r_{n})$. It follows that
$$|\{1\leq i\leq n-1 : (r_{i}-r_{n}, k)\neq1\}|\leq1.$$
 Since $|\hat{X_{n}}|<p_{1}$, in view of
Lemma \ref{bing}, we have
$$|(\hat{X_{n}}+(r_{s}-r_{n}))\setminus\hat{X_{n}}|\geq1.$$
Then using Lemma \ref{gaijinchowla}, we get
$$|\hat{X_{n}}+(\hat{A_{1}}\cup\ldots\cup\hat{A_{t}}\cup\hat{A_{n}})|
\geq |\hat{X_{n}}|+t \quad\text{for}\quad 1\leq t\leq n-1,$$
and consequently,
\begin{align*}
  |\triangle _{nn}|&=|(X_{n}+ A )\backslash(X_{n}+ A_{n})|\geq|A_{1}|+\ldots+|A_{n-1}|.
\end{align*}

Combining the above we have completed the proof.
\end{proof}

\begin{lemma}
\label{2pxiaoyinli}  Let $A$ be a nonempty finite subset of $\mathbb{Z}$.
Then
$$|A+k\cdot A|\geq(k+1)|A|-k!.$$
\end{lemma}

\begin{proof} We use induction to show that
\begin{equation}\label{t chenli}
|A+k\cdot A|\geq(t+1)|A|-(t-1)!k.
\end{equation}
holds for every $t=1,\ldots,k$.

Clearly (\ref{t chenli}) is true for $t=1$.

Now assume that (\ref{t chenli}) holds for a fixed $1\le t<k$. We want to deduce (\ref{t chenli}) with $t$ replaced by $t+1$, i.e.,
\begin{equation}\label{t2}
|A+k\cdot A|\geq (t+2)|A|-t!k.
\end{equation}
As discussed in Lemma \ref{small}, we only need to deal with the case $h\leq t$. By Lemma \ref{a+kb} and the induction hypothesis, we have
\begin{equation}\label{ff}
|A_i+k\cdot A|\geq |X_i+A_{1}|\geq (t+1)|A_i|+|A_1|-k\quad\text{for}\quad
i\in F
\end{equation}
and
\begin{equation}\label{ee}
|A_i+k\cdot A|=|X_i+A|\geq|X_i+k\cdot X_{i}|+|\triangle_{ii}|\geq
(t+1)|A_i|-(t-1)!k+|\triangle_{ii}|\quad\text{for}\quad i\in E.
\end{equation}

{\it Case} 1. $(r_2, k)=1$.

If $2\in F$, then $2\not\in E$. By Lemma \ref{changsui}, we have
$|\triangle_{ii}|\geq|A_2|$ for $i\in E$. Combining (\ref{ff}) and
(\ref{ee}), we have
\begin{align*}
|A+k\cdot A|&=\sum\limits_{i\in F}|A_{i}+k\cdot A|+\sum\limits_{i\in
E}|A_{i}+k\cdot A|\\
&\geq (t+1)|A|+|F||A_1|+|E||A_2|-\(k|F|+|E|(t-1)!k\)\\
&\geq (t+2)|A|-t!k.
\end{align*}

When $2\in E$, we have $|\triangle_{22}|\geq|A_1|$. Furthermore,
$|\triangle_{ii}|\geq|A_2|$ for $i\in E\setminus\{2\}$. Hence
\begin{align*}
|A+k\cdot A|&=\sum\limits_{i\in F}|A_{i}+k\cdot A|+\sum\limits_{i\in
E}|A_{i}+k\cdot A|\\
&\geq (t+1)|A|+|F||A_1|+|A_1|+(|E|-1)|A_2|-\(k|F|+|E|(t-1)!k\)\\
&\geq (t+2)|A|-t!k.
\end{align*}

{\it Case} 2. $(r_2, k)=p_1$.

 Observe that
\begin{align*}
|A+k\cdot A|&=|A_n+k\cdot A|+\sum\limits_{i\in
F\setminus\{n\}}|A_{i}+k\cdot A|+\sum\limits_{i\in
E\setminus\{n\}}|A_{i}+k\cdot A|\\
&\geq |X_n+A|+\sum\limits_{i\in
F\setminus\{n\}}(|A_{i}|+k|A_{1}|-k)+\sum\limits_{i\in E\setminus\{n\}}|A_{i}+k\cdot A_{i}|+\sum\limits_{i\in
E\setminus\{n\}}|\triangle_{ii}|\\
&\geq |X_n+A|+(t+1)\sum\limits_{i\neq
n}|A_i|+|F\setminus\{n\}||A_1|+\sum\limits_{i\in
E\setminus\{n\}}|\triangle_{ii}|\\
&\ \ \ \ -\(k|F\setminus\{n\}|+|E\setminus\{n\}|(t-1)!k\)\\
&\geq |X_n+A|+|F\setminus\{n\}||A_1|+\sum\limits_{i\in
E\setminus\{n\}}|\triangle_{ii}|+(t+1)\sum\limits_{i\neq
n}|A_i|-(t-1)(t-1)!k.
\end{align*}
For convenience we denote
$$\mathcal{S}=|X_n+A|+|F\setminus\{n\}||A_1|+\sum\limits_{i\in
E\setminus\{n\}}|\triangle_{ii}|.$$ In order to get (\ref{t2}), it
is sufficient to prove
$$\mathcal{S}\geq (t+1)|A_n|+|A|-(t-1)!k.$$

If $n\in F$, then
$$|X_{n}+A|\geq |X_n|+k|A_1|-k\geq
(t+1)|A_n|+t(|A_1|-|A_n|)+|A_1|-k.$$ Notice that
$|\triangle_{ii}|\geq|A_n|$ for every $i\in E$ and that $h\leq t$. Thus
\begin{align*}
\mathcal{S}&\geq(t+1)|A_n|+t(|A_1|-|A_n|)+|F||A_1|+|E||A_n|-k\\
&\geq (t+1)|A_n|+|A|-(t-1)!k.
\end{align*}

Below we assume $n\in E$.

When $|A_1|\leq p_2|A_n|$, by Lemma \ref{changsui} we
get
\begin{align*}
\mathcal{S}&\geq
(t+1)|A_n|-(t-1)!k+|A_2|+|F||A_1|\\
&\ \ \ \ +|A_2||E\cap\{1\}|+|A_1||E\cap\{2,\cdots,n-1\}|+|A_n||E\cap\{n+1,\cdots,h\}|\\
&\geq (t+1)|A_n|+|A|-(t-1)!k.
\end{align*}

Now suppose $|A_1|> p_2|A_n|$. If $|\hat{X_n}|<p_1$, then from Lemma \ref{changsui}
\begin{align*}
\mathcal{S}&\geq
(t+1)|A_n|-(t-1)!k+|A_1|+\cdots+|A_{n-1}|+|F||A_1|+(|E|-1)|A_n|\\
&\geq (t+1)|A_n|+|A|-(t-1)!k.
\end{align*}

If $|\hat{X_n}|\geq p_1>p_2$, then
\begin{align*}
|X_n+A|&\geq
|A_n|+p_1|A_1|-k=(k+1)|A_n|+p_1(|A_1|-p_2|A_n|)-k\\
&\geq(t+1)|A_n|+|A_n|+(n-1)(|A_1|-p_2|A_n|)-k
\end{align*}
since $p_1> p_2\geq n-1$. When $|\hat{X_n}|\geq p_2>p_1$ we also have
\begin{align*}
|X_n+A|&\geq
|A_n|+p_2|A_1|-k=(k+1)|A_n|+p_2(|A_1|-p_1|A_n|)-k\\
&\geq(t+1)|A_n|+|A_n|+(n-1)(|A_1|-p_2|A_n|)-k.
\end{align*}
With the help of Lemma \ref{changsui}, we deduce that
$$\mathcal{S}\geq (t+1)|A_n|+|A|-(t-1)!k$$ when $|\hat{X_n}|\geq
p_1>p_2$ or $|\hat{X_n}|\geq p_2>p_1$.

If $p_1\leq |\hat{X_n}|<p_2$, then
\begin{align*}
\mathcal{S}&\geq
|A_n|+|\hat{X_n}||A_1|+\sum_{i=2}^{l}|A_i|-k+|F||A_1|\\
&\ \ \ \ +p_2|A_n||E\cap\{2,\cdots,n-1\}|+|A_n||E\cap\{n+1,\cdots,h\}|\\
&\geq
\sum_{i=1}^{n}|A_i|+(|\hat{X_n}|-n+l)|A_1|+(n-2)p_2|A_n|+(h-n)|A_n|-k\\
&\geq |A|+(|\hat{X_n}|+l-2)p_2|A_n|-k.
\end{align*}
Observe that
$$|\hat{X_n}|+l-2=\min\{|\hat{X_n}|+n-3, p_2-1\}\geq
p_1.$$ Therefore
$$\mathcal{S}\geq |A|+k|A_n|-k\geq
|A|+(t+1)|A_n|-(t-1)!k.$$

This concludes the proof.
\end{proof}

\begin{proof}[Proof of Theorem~\ref{liangsushuchenji}]
Suppose $|A|\geq (k-1)^{2}k!$. As discussed in the proof of
Theorem \ref{t1}, (\ref{2}) is valid when $h=k$ or $|\hat{X_{i}}|=k$
for all $1\leq i\leq h$. Below assume $h<k$ and
$E\neq\emptyset$. Then $|A_{1}|\geq |A|/h\geq hk!$.

{\it Case} 1.  $(r_{2},k)=1$.

If $|\hat{X_{2}}|=k$,  then there exists $s\in E\setminus\{2\}$ since $E\neq\emptyset$.
From Lemma \ref{changsui},
$|\triangle_{ss}|\geq|A_{2}|$. Then
in light of Lemmas \ref{a+kb} and \ref{2pxiaoyinli}
\begin{align*}
|A+k\cdot A|&=|A_{2}+k\cdot A|+|A_{s}+k\cdot A|+\sum\limits_{i\neq 2,s}|A_{i}+k\cdot A|\\
&\geq|A_{2}|+k|A_{1}|-k+(k+1)|A_{s}|-k!+|A_{2}|+\sum\limits_{i\neq 2,s}\((k+1)|A_{i}|-k!\)\\
&\geq(k+1)|A|-hk!+|A_{1}|\geq(k+1)|A|.
\end{align*}

If $|\hat{X_{2}}|<k$, then $|\triangle_{22}|\geq|A_{1}|$ and
\begin{align*}
|A+k\cdot A|&=|A_{2}+k\cdot A|+\sum\limits_{i\neq 2}|A_{i}+k\cdot A|\\
&=|A_{2}+k\cdot A_{2}|+|\triangle_{22}|+\sum\limits_{i\neq 2}|A_{i}+k\cdot A|\\
&\geq(k+1)|A_{2}|-k!+|A_{1}|+\sum\limits_{i\neq 2}\((k+1)|A_{i}|-k!\)\\
&\geq(k+1)|A|-hk!+|A_{1}|\geq(k+1)|A|.
\end{align*}

{\it Case} 2. $(r_{2},k)=p_1$.

When $|\hat{X_{n}}|=k$, using Lemma \ref{changsui} we have $s\in E$ with $s\neq n$ such that
 $|\triangle_{ss}|\geq |A_{n}|$, which states
\begin{align*}
|A+k\cdot A|&=|A_{n}+k\cdot A|+|A_{s}+k\cdot A|+\sum\limits_{i\neq ,n,s}|A_{i}+k\cdot A|\\
&\geq|A_{n}|+k|A_{1}|-k+(k+1)|A_{s}|-k!+|A_{n}|+\sum\limits_{i\neq n,s}\((k+1)|A_{i}|-k!\)\\
&\geq(k+1)|A|.
\end{align*}

Below suppose $|\hat{X_{n}}|<k$.

{\it Subcase} 1. $|\hat{X_{n}}|\geq p_{1}$ and $ |A_{1}|\leq p_{2}|A_{n}|$.

In the case $|\hat{X_{2}}|=k$, by Lemma \ref{changsui} we have $|\triangle_{nn}|\geq
|A_{2}|$ and
\begin{align*}
|A+k\cdot A|&=|A_{2}+k\cdot A|+|A_{n}+k\cdot A|+\sum\limits_{i\neq 2,n}|A_{i}+k\cdot A|\\
&\geq|A_{2}|+k|A_{1}|-k+(k+1)|A_{n}|-k!+|A_{2}|+\sum\limits_{i\neq2,n}\((k+1)|A_{i}|-k!\)\\
&\geq(k+1)|A|.
\end{align*}

When $|\hat{X_{2}}|<k$, we obtain $|\triangle_{22}|\geq|A_{1}|$ from Lemma
\ref{changsui} and hence
\begin{align*}
|A+k\cdot A|&=|A_{2}+k\cdot A_{2}|+|\triangle_{22}|+\sum\limits_{i\neq 2}|A_{i}+k\cdot A|\\
&\geq(k+1)|A_{2}|-k!+|A_{1}|+\sum\limits_{i\neq 2}\((k+1)|A_{i}|-k!\)\\
&\geq(k+1)|A|.
\end{align*}

{\it Subcase} 2. $|\hat{X_{n}}|\geq p_{1}$ and $ |A_{1}|> p_{2}|A_{n}|$.

If $|\hat{X_{n}}|> p_{1}$, then
$$|X_{n}+A|\geq|A_{n}|+|\hat{X_{n}}||A_{1}|-k\geq|A_{n}|+p_{1}|A_{1}|+|A_{2}|-k.$$
For $|\hat{X_{n}}|=p_{1}$, using Lemma \ref{bing}, we have
$|(\hat{X_{n}}+\hat{A_{2}})\backslash \hat{X_{n}}|\geq 1$  and then
$$|X_{n}+A|\geq|X_{n}+A_{1}|+|A_{2}|\geq|A_{n}|+p_{1}|A_{1}|+|A_{2}|-k.$$

Applying the above, if $|\hat{X_{2}}|=k$, then
\begin{align*}
|A+k\cdot A|&=|A_{2}+k\cdot A|+|A_{n}+k\cdot A|+\sum\limits_{i\neq 2,n}|A_{i}+k\cdot A|\\
&\geq|A_{2}|+k|A_{1}|-k+|A_{n}|+p_{1}|A_{1}|+|A_{2}|-k+\sum\limits_{i\neq 2, n}\((k+1)|A_{i}|-k!\)\\
&\geq(k+1)|A|.
\end{align*}

In the case $|\hat{X_{2}}|<k$, clearly $|\triangle_{22}|\geq p_2|A_{n}|$ from Lemma \ref{changsui} and therefore,
\begin{align*}
|A+k\cdot A|&=|A_{2}+k\cdot A|+|A_{n}+k\cdot A|+\sum\limits_{i\neq 2,n}|A_{i}+k\cdot A|\\
&\geq|\triangle_{22}|+|A_{n}|+p_{1}|A_{1}|+|A_{2}|-k+\sum\limits_{i\neq n}\((k+1)|A_{i}|-k!\)\\
&\geq(k+1)|A|.
\end{align*}

{\it Subcase} 3. $|\hat{X_{n}}|< p_{1}$.

In this case we get
\begin{align*}
|A+k\cdot A|&=|A_{n}+k\cdot A_{n}|+|\triangle_{nn}|+\sum\limits_{i\neq n}|A_{i}+k\cdot A|\\
&\geq|\triangle_{nn}|+\sum\limits_{1\leq i\leq
h}\((k+1)|A_{i}|-k!\)\geq(k+1)|A|
\end{align*}
because of $|\triangle_{nn}|\geq |A_{1}|$ from Lemma \ref{changsui}.

Combining the above we complete the proof of Theorem \ref{liangsushuchenji}.
\end{proof}

\section{Proof of Theorem \ref{k=4}}

\begin{lemma}
\label{234} Let $A\subseteq \mathbb{Z}$, then\\
(i) \ \ $|A+4\cdot A|= 4$, if $|A|=2$. \\
(ii)\ \ $|A+4\cdot A|\geq 8$, if $|A|=3$. \\
(iii)\ $|A+4\cdot A|\geq 12$, if $|A|=4$.
\end{lemma}
\begin{proof}
Lemma \ref{234} can be proved easily by a direct analysis.
\end{proof}

Observing that $$|(X_{i}+A_{1})\setminus(X_{i}+A_{2})|\geq1 \quad\text{or}\quad
|(X_{i}+A_{2})\setminus(X_{i}+A_{1})|\geq1$$
when $|A_{1}|=|A_{2}|$,
so in this case, we may suppose $|(X_{i}+A_{2})\setminus(X_{i}+A_{1})|\geq1$ without loss of generality.
Below we fix $A\subseteq \mathbb{Z}$ with $|A|\geq 5$, and use the notations in Section 2.
\begin{lemma}
\label{x3} When $h=3$, for $i=1,2,3$ we have $$|A_{i}|\leq 4 \Rightarrow |A_{i}+4\cdot A|\geq |A|+2|A_{i}|-2.$$
\end{lemma}
\begin{proof}
Recall that $A_{i}=4\cdot X_{i}+r_{i}$,
$|A_{1}|\geq|A_{2}|\geq|A_{3}|$ and $|A_{i}+4\cdot A|=|X_{i}+A|$.

(I) If $|\hat{X_{i}}|=1$, in light of Lemma \ref{a+kb} we have
\begin{align*}
  |X_{i}+A|&=|X_{i}+A_{1}|+|X_{i}+A_{2}|+|X_{i}+A_{3}|\\
  &\geq|X_{i}|+|A_{1}|-1+|X_{i}|+|A_{2}|-1+|X_{i}|+|A_{3}|-1\\
  &\geq |A|+3|A_{i}|-3\\
  &\geq |A|+2|A_{i}|-2.
\end{align*}

(II) When $|\hat{X_{i}}|=2$, Lemma \ref{chowla} implies
$|\hat{X_{i}}+\hat{A}|=4$.

In the case $|\hat{X_i}+\hat{A_1}\cup\hat{A_2}|\geq 3$,
\begin{align*}
|X_i+A|&=|X_i+A_1|+|(X_i+A)\setminus(X_i+A_1)|\\
&\geq |A_i|+2|A_1|-2+|A_2|+|A_3|\\
&\geq |A|+2|A_i|-2.
\end{align*}

When $|\hat{X_i}+\hat{A_1}\cup\hat{A_2}|=2$, we have
\begin{align*}
|X_i+A|&=|X_i+A_1\cup A_2|+|X_i+A_3|\\
&=|X_i+A_1|+|(X_i+A_2)\setminus(X_i+A_1)|+|X_i+A_3|\\
&\geq |A_i|+2|A_1|-2+\delta_{|A_{1}|,\ |A_{2}|}+|A_i|+2|A_3|-2\\
&\geq |A|+2|A_i|-2,
\end{align*}
where the Kronecker symbol $\delta_{s,t}$ takes $1$ or $0$ according to $s=t$ or not.

(III) If $|\hat{X_{i}}|=3$, then
$$|X_{i}+A_{1}|\geq 3|A_{1}|\ \text{for}\ |A_{i}|=3\quad \text{and}\quad
|X_{i}+A_{1}|\geq 3|A_{1}|+1\ \text{for}\ |A_{i}|=4.$$

Clearly $|(\hat{X_i}+\hat{A_2})\setminus(\hat{X_i}+\hat{A_1})|\geq 1$ and hence in the case $|A_{1}|>|A_{3}|$ we have
$$|X_{i}+A|\geq
|X_{i}+A_{1}|+|A_{2}|\geq |A|+2|A_{i}|-2.$$

In the case
$|A_{1}|=|A_{2}|=|A_{3}|$, the reader can get directly that
$$|X_{i}+A|\geq 13\ \text{if}\ |A_{i}|=3\quad \text{and}\quad
|X_{i}+A|\geq 18\ \text{if}\ |A_{i}|=4.$$

Thus we also have
$|X_{i}+A|\geq |A|+2|A_{i}|-2$.

(IV)
 If $|\hat{X_{i}}|=4$, then we have
$|A_{i}|=4$ since $|A_{i}|\leq 4$. Let $a$ be the minimal number of $X_{i}$ and let $b$ be
the maximal one of $A$. Because $|\hat{X_{i}}|=4$ and
$|\hat{A}|=3$, we have $|\hat{\{a\}}+\hat{A}|=3$ and
$|\hat{X_{i}}+(\widehat{A_{i}\setminus \{b\}})|=4$. Then
$$|(X_{i}+(A_{i}\backslash\{b\}))\backslash((a+A)\cup(b+X_{i}))|\geq3.$$ It turns out
that
\begin{align*}
|X_{i}+A|&\geq|a+A|+|b+X_{i}|-1+|(X_{i}+(A_{i}\backslash\{b\}))\backslash((a+A)\cup(b+X_{i}))|\\
&\geq|A|+4-1+3\geq |A|+2|A_{i}|-2.
\end{align*}

The proof is complete.
\end{proof}

\begin{lemma}
\label{x2}When $h=2$, we have for $i=1, 2$ that
$$|A_{i}|\leq 4\Rightarrow|A_{i}+4\cdot A|\geq |A|+3|A_{i}|-3.$$
\end{lemma}

\begin{proof}
We divide the proof into four parts.

(I) When $|\hat{X_{i}}|=1$, we use Lemma \ref{234} to obtain
for $|A_{i}|\leq 4$ that,
\begin{align*}
|X_{i}+A|&=|X_{i}+A_{i}|+|X_{i}+(A\backslash A_{i})|\\
&\geq |X_{i}+4\cdot X_{i}|+|A|-|A_{i}|+|A_{i}|-1\\
&\geq |A|+3|A_{i}|-3.
\end{align*}

(II) If $|\hat{X_{i}}|=2$, then we have $3
\leq |\hat{X_{i}}+\hat{A}| \leq4$ from Lemma \ref{chowla}. Then we distinguish two
cases.

{\it Case} 1. $|\hat{X_{i}}+\hat{A}|= 4$.

Since $|\hat{X_{i}}+\hat{A_{i}}|= 2$, by Lemma \ref{234} and
$|A|\geq 5$, we get
\begin{align*}
|X_{i}+A|&=|X_{i}+A_{i}|+|X_{i}+(A\backslash A_{i})|\\
&\geq |X_{i}+4\cdot X_{i}|+2(|A|-|A_{i}|)\\
&\geq |A|+3|A_{i}|-3.
\end{align*}

{\it Case} 2. $|\hat{X_{i}}+\hat{A}|= 3$.

Define $\eta=|\  \{x\in X_{i}:(\hat{\{x\}}+\widehat{A \backslash
A_{i}})\nsubseteqq (\hat{X_{i}}+\hat{A_{i}}) \}\ |$. When
$\eta\geq2$, by Lemma \ref{234} we have
\begin{align*}
|X_{i}+A|&\geq|X_{i}+A_{i}|+|A|-|A_{i}|+1\\
&\geq |A|+3|A_{i}|-3.
\end{align*}

Below suppose
 $\eta=1$. It is easy to see that
$$|X_{i}+A|\geq|A|+|A|-|A_{i}|\geq|A|+3|A_{i}|-3$$ for $|A_{i}|=2$ since $|A|\geq 5$.

For $|A_{i}|=3$, we write $X_{i}=\{s_{1},s_{2},s_{3}\}$ with
$s_{1}\equiv s_{2}\ \ (\rm mod\  \ 4 )$ and $s_{1}<s_{2}$. Then we have
$\hat{\{s_{3}\}}+\hat{A_{i}}=\hat{\{s_{1}\}}+\hat{A}\backslash\hat{A_{i}}$ because of $\eta=1$. Now we show
\begin{equation}\label{duo yigeyusu dingli}
|((s_{1}\cup s_{2})+A\backslash
A_{i})\backslash(s_{3}+A_{i})|\geq1.
\end{equation}
Clearly (\ref{duo yigeyusu dingli})
holds for $|A|\geq 6$. If $|A|=5$, then $|A_{1}|=3$. Let
$A_{1}=\{a_{1}, a_{2},a_{3}\}$ with $a_{1}=4s_{1},\  a_{2}=4s_{2}$
and $\ a_{3}=4s_{3}$. And let $A_{2}=\{a_{4},a_{5}\}$ with $a_{4}<a_{5}$.
If $a_{1}<a_{2}<a_{3}$ and (\ref{duo yigeyusu dingli}) fails, then
$s_{3}+a_{1}=s_{1}+a_{4}, \ s_{3}+a_{2}=s_{1}+a_{5}=s_{2}+a_{4}$
and $s_{3}+a_{3}=s_{2}+a_{5}$, and hence $a_{2}-a_{1}=s_{2}-s_{1}$, which contradicts
$a_{1}-a_{2}=4(s_{1}-s_{2})$. When
$a_{3}<a_{1}<a_{2}$ or $a_{1}<a_{3}<a_{2}$, we also get
(\ref{duo yigeyusu dingli}). From (\ref{duo yigeyusu dingli}) and Lemma \ref{234} we
have
\begin{align*}
|X_{i}+A|&\geq|X_{i}+A_{i}|+|A|-|A_{i}|+1\\
&\geq 8+|A|-3+1\geq |A|+3|A_{i}|-3.
\end{align*}

For $|A_{i}|=4$, we have
$|X_{i}+4\cdot X_{i}|\geq 13$ in the case $|\hat{X_{i}}|=2$ and $\eta=1$, and therefore
\begin{align*}
|X_{i}+A|&\geq |X_{i}+A_{i}|+|A|-|A_{i}|\\
&\geq 13+|A|-4\geq |A|+3|A_{i}|-3.
\end{align*}
Now we give the reason for $|X_{i}+4\cdot X_{i}|\geq 13$. We write $X_{i}=X_{i1}\cup X_{i2}$
with $X_{i1}=4\cdot Y_{1}+r_{1}$ and $X_{i2}=4\cdot Y_{2}+r_{2}$. The fact $\eta=1$ allows us to
assume $|Y_{1}|=3$ and $|Y_{2}|=1$. To discuss $|X_{i}+4\cdot X_{i}|$, we may suppose $0\in \hat{X_{i}}$ and
$\gcd (X_{i})=1$ without loss of generality. Then with the help of Lemma \ref{chowla}
we have $|(\hat{Y_{1}}+\hat{X_{i2}})\backslash(\hat{Y_{1}}+\hat{X_{i1}})|\geq 1$ since $|\hat{Y_{1}}|\leq 3$ and hence
\begin{align*}
|X_{i}+4\cdot X_{i}|&=
|Y_{1}+X_{i}|+|Y_{2}+X_{i}|\\
&\geq|Y_{1}+X_{i1}|+1+|X_{i}|\\
&=|Y_{1}+4\cdot Y_{1}|+1+|X_{i}|\\&\geq 8+1+4=13.
\end{align*}

(III) In the case $|\hat{X_{i}}|$=3, we have $|\hat{X_{i}}+\hat{A}|=4$.
Applying Lemma \ref{234} we get
\begin{align*}
|X_{i}+A|&\geq
|X_{i}+A_{i}|+|A|-|A_{i}|\\
&\geq3|A_{i}|+|A_{i}|-3+|A|-|A_{i}|\geq |A|+3|A_{i}|-3.
\end{align*}

(IV) If $|\hat{X_{i}}|=4$, then $|A_{i}|=4$. For $|A|=8$,
we have $|A_{1}|=|A_{2}|=4$ and
 $|(X_{i}+A_{2})\backslash(X_{i}+A_{1})|\geq 1$ by the assumption.
Note that $ |X_{i}+A_{1}|=4|A_{1}|$ since $|\hat{X_{i}}|=4$. Thus
\begin{align*}
|X_{i}+A|&\geq|X_{i}+A_{1}|+|(X_{i}+A_{2})\backslash(X_{i}+A_{1})|\\
&\geq 4|A_{1}|+\delta_{|A_{1}|, |A_{2}|}\\
&\geq|A|+3|A_{i}|-3.
\end{align*}
\end{proof}

\begin{proof}[Proof of Theorem~\ref{k=4}]
If $|\hat{A}|=4$, then $|A+4\cdot A|\geq 5|A|-4$ and hence (\ref{3}) holds. Below we
 assume $|\hat{A}|\leq 3$.

We prove (\ref{3}) by induction on $|A|$. Clearly, (\ref{3}) holds for $|A|=5$ with the help of Lemmas \ref{x3} and \ref{x2}.
Now we let $|A|>5$ and assume that
$$|B+4\cdot B|\geq 5|B|-6\ \ \text{for any}\ \  B\subset\Z\ \ \text{with}\ \  5 \leq |B|<|A|.$$
We divide our proof of (\ref{3}) into three parts.

{\bf Claim I}. (\ref{3}) holds when $h=3$ and $|A_{3}|\leq4$.

To prove Claim I, we distinguish three small cases.

{\it Case} I.1. $|A_{1}|\leq 4$.

We use Lemma \ref{x3} to obtain
$$|A+4\cdot A|=|X_{1}+A|+|X_{2}+A|+|X_{3}+A|\geq 5|A|-6.$$

{\it Case} I.2. $|A_{1}|\geq 5$ and $|A_{2}|\leq 4$.

For $|\hat{X_{1}}|=4$, by Lemmas \ref{a+kb} and \ref{x3}, we have
\begin{align*}
|A+4\cdot A|&=|X_{1}+A|+|X_{2}+A|+|X_{3}+A|\\
&\geq|X_{1}+A_{1}|+|X_{2}+A|+|X_{3}+A|\\
&\geq5|A_{1}|-4+|A|+2|A_{2}|-2+|A|+2|A_{3}|-2\\
&\geq5|A|+|A_{1}|-|A_{2}|+|A_{1}|-|A_{3}|-8\geq 5|A|-6.
\end{align*}

For $|\hat{X_{1}}|\leq 3$, since $|\hat{A}|=3$, applying Lemma
\ref{chowla} we have $|\hat{X_{1}}+\hat{A}|>|\hat{X_{1}}|$ and
then
\begin{align*}
|A+4\cdot A|&=|X_{1}+A|+|X_{2}+A|+|X_{3}+A|\\
&\geq|X_{1}+4 \cdot X_{1}|+|A_{3}|+|X_{2}+A|+|X_{3}+A|\\
&\geq5|A_{1}|-6+|A_{3}|+|A|+2|A_{2}|-2+|A|+2|A_{3}|-2\\
&\geq5|A|+|A_{1}|-|A_{2}|+|A_{1}|-10\geq 5|A|-6.
\end{align*}

{\it Case} I.3. $|A_{2}|\geq 5$ and $|A_{3}|\leq 4$.

We have $|\triangle_{11}|+|\triangle_{22}|\geq 2$ by Lemma \ref{ii} and hence
\begin{align*}
|A+4\cdot A|&=|X_{1}+A|+|X_{2}+A|+|X_{3}+A|\\
&\geq|X_{1}+4 \cdot X_{1}|+|\triangle_{11}|+|X_{2}+4 \cdot
X_{2}|+|\triangle_{22}|+|X_{3}+A|\\
&\geq5|A_{1}|-4+5|A_{2}|-4+2+|A|+2|A_{3}|-2\\
&\geq5|A|+|A_{1}|-|A_{3}|+|A_{2}|-|A_{3}|-8\geq 5|A|-6
\end{align*}
when $|\hat{X_{1}}|=4$ and $|\hat{X_{2}}|=4$.
In the case $|\{1\le i\le 2 :\, |\hat{X_{i}}|=4\}|=1$ we obtain
\begin{align*}
|A+4\cdot A|&=|X_{1}+A|+|X_{2}+A|+|X_{3}+A|\\
&\geq5(|A_{1}|+|A_{2}|)-6-4+|A_{3}|+|A|+2|A_{3}|-2\\
&\geq5|A|+|A_{1}|+|A_{2}|-|A_{3}|-12\geq 5|A|-6
\end{align*}
in view of
Lemma \ref{chowla} and the induction hypothesis.

When $|\hat{X_{1}}|\leq3$ and $|\hat{X_{2}}|\leq3$, by Lemma
\ref{chowla} we have $|\hat{X_{1}}+\hat{A}|>|\hat{X_{1}}|$ and
$|\hat{X_{2}}+\hat{A}|>|\hat{X_{2}}|$. Then
\begin{align*}
|A+4\cdot A|&=|X_{1}+A|+|X_{2}+A|+|X_{3}+A|\\
&\geq5|A_{1}|-6+|A_{3}|+5|A_{2}|-6+|A_{3}|+|A|+2|A_{3}|-2\\
&\geq5|A|+|A_{1}|+|A_{2}|-14\geq 5|A|-6.
\end{align*}

{\bf Claim II.} (\ref{3}) holds when $h=3$ and $|A_{3}|\geq 5$.

By Lemma \ref{chowla}, when $|\hat{X_{i}}|\leq 3$, we have $|\hat{X_{i}}+\hat{A}|>|\hat{X_{i}}|$ and hence
$$|X_{i}+A|\geq |X_{i}+A_{i}|+|A_{3}|.$$

If $|\hat{X_{1}}|\leq 3$, $|\hat{X_{2}}|\leq 3$ and
$|\hat{X_{3}}|\leq 3$, then we have
\begin{align*}
|A+4\cdot A|&=|X_{1}+A|+|X_{2}+A|+|X_{3}+A|\\
&\geq5|A_{1}|-6+|A_{3}|+5|A_{2}|-6+|A_{3}|+5|A_{3}|-6+|A_{3}|\\
&\geq5|A|+3|A_{3}|-18\geq 5|A|-6.
\end{align*}

When $|\ \{i:|\hat{X_{i}}|\leq 3\}\ |=2$, by Lemma
 \ref{a+kb} we get
$$|A+4\cdot A|\geq 5(|X_{1}|+|X_{2}|+|X_{3}|)-6-6-4+2|A_{3}|\geq 5|A|-6.$$

In the case $|\hat{X_{1}}|=|\hat{X_{2}}|=|\hat{X_{3}}|=4$, we have
$|\triangle_{11}|+|\triangle_{22}|+|\triangle_{33}|\geq 6$ in light
of Lemma \ref{ii}. Then by
Lemma \ref{a+kb} we get
\begin{align*}
|A+4\cdot A|&=|X_{1}+A|+|X_{2}+A|+|X_{3}+A|\\
&\geq|X_{1}+A_{1}|+|\triangle_{11}|+|X_{2}+A_{2}|+|\triangle_{22}|+|X_{3}+A_{3}|+|\triangle_{33}|\\
&\geq5(|A_{1}|+|A_{2}|+|A_{3}|)-4-4-4+6\geq 5|A|-6.
\end{align*}

It remains to handle the case $|\ \{i:|\hat{X_{i}}|\leq 3\}\ |=1$,
and we make detailed discussions.

{\it Case} II.1. $|A_{1}|=|A_{2}|=|A_{3}|$.

We may suppose $|\hat{X_{3}}|\leq 3$ and
$|\hat{X_{1}}|=|\hat{X_{2}}|=4$. Note that $|X_{1}+A_{i}|\geq5|A_{1}| -4$ and $|X_{2}+A_{i}|\geq5|A_{2}| -4$ for all $i$. Now we prove
$$|X_{1}+A|\geq5|A_{1}| -2\ \ \text{and}\ \  |X_{2}+A|\geq5|A_{2}| -2.$$ Note that
 $$\min (X_{1}+A_{i})\notin X_{1}+A_{j}\ \ \text{or}\ \  \max (X_{1}+A_{i})\notin
 X_{1}+A_{j}$$
 if $|X_{1}+A_{i}|=|X_{1}+A_{j}|$. When $|X_{1}+A_{i}|\geq 5|A_{1}|-3$ for some $i$,
 we have $|X_{1}+A|\geq5|A_{1}|-2$ since $|X_{1}+A\setminus A_{i}|\geq5|A_{1}|-3$. If $|X_{1}+A_{i}|= 5|A_{1}|-4$ for all $i$,
 then $A_{1}, A_{2}$ and $A_{3}$
  must be arithmetic progressions with the same difference by Lemma \ref{a+kb}, and therefore $|\triangle_{11}|\geq 2$
  and $$|X_{1}+A|\geq|X_{1}+A_{1}|+|\triangle_{11}|\geq5|A_{1}|-2.$$ Similarly,
$|X_{2}+A|\geq5|A_{2}|-2$. So
\begin{align*}
|A+4\cdot A|&=|X_{1}+A|+|X_{2}+A|+|X_{3}+A|\\
&\geq5|A_{1}|-2+5|A_{2}|-2+5|A_{3}|-6+|A_{3}|\\
&\geq5|A|+|A_{3}|-10\\
&\geq 5|A|-6.
\end{align*}

{\it Case} II.2. $|A_{1}|=|A_{2}|=|A_{3}|$ fails.

Note that $|A_{1}|>|A_{3}|\geq 5$. If $|A_{1}|>|A_{2}|$, then
$$|X_{2}+A|\geq
|X_{2}+4\cdot X_{1}|\geq|A_{2}|+4|A_{1}|-4\geq5|A_{2}|$$
or
$$|X_{3}+A|\geq |X_{3}+4\cdot
X_{1}|\geq|A_{3}|+4|A_{1}|-4\geq5|A_{3}|$$
since
$|\ \{i:|\hat{X_{i}}|\leq 3\}\ |=1$ . Hence
\begin{align*}
|A+4\cdot A|&=|X_{1}+A|+|X_{2}+A|+|X_{3}+A|\\
&\geq5(|A_{1}|+|A_{2}|+|A_{3}|)-6-4+|A_{3}|\\
&\geq 5|A|-6.
\end{align*}

Now suppose $|A_{1}|=|A_{2}|>|A_{3}|$. If $|\hat{X_{3}}|\leq 3$, then $|\hat{X_{1}}|=|\hat{X_{2}}|=4$. As mentioned in case II.1 we have
$$|A+4\cdot A|\geq
5|A_{1}|-4+1+5|A_{2}|-4+1+5|A_{3}|-6+|A_{1}|\geq 5|A|-6.$$ When
$|\hat{X_{3}}|=4$, it is clear that $$|A+4\cdot A|\geq
5(|A_{1}|+|A_{2}|)-4-6+|A_{3}|+|A_{3}|+4|A_{1}|-4\geq 5|A|-6.$$

{\bf Claim III.} (\ref{3}) holds for $h=2$.

We first note that if $|\hat{X_{i}}|\leq 3$ then
$|\hat{X_{i}}+\hat{A}|>|\hat{X_{i}}+\hat{A_{i}}|$ by Lemma \ref{chowla}.

{\it Case} III.1. $|A_{1}|\leq 4$.

Applying Lemma \ref{x2} we get
\begin{align*}
|A+4\cdot A|&=|X_{1}+A|+|X_{2}+A|
\\
&\geq|A|+3|A_{1}|-3+|A|+3|A_{2}|-3 \geq 5|A|-6.
\end{align*}

{\it Case} III.2. $|A_{1}|\geq 5$ and $|A_{2}|\leq 4$.

When $|\hat{X_{1}}|=4$, we have
\begin{align*}
|A+4\cdot A|&=|X_{1}+A|+|X_{2}+A|\geq|X_{1}+A_{1}|+|X_{2}+A|\\
&\geq5|A_{1}|-4+|A|+3|A_{2}|-3\\
&=5|A|+|A_{1}|-|A_{2}|-7
\geq5|A|-6.
\end{align*}

If $|\hat{X_{1}}|\leq 3$, then
\begin{align*}
|A+4\cdot
A|&=|X_{1}+A|+|X_{2}+A|\\
&\geq|X_{1}+A_{1}|+|\triangle_{11}|+|X_{2}+A|\\
&\geq5|A_{1}|-6+|A_{2}|+|A|+3|A_{2}|-3
\\&=5|A|+|A_{1}|-|A_{2}|+|A_{2}|-9 \geq5|A|-6.
\end{align*}

{\it Case} III.3. $|A_{2}|\geq 5$.

For $|\hat{X_{1}}|=|\hat{X_{2}}|=4$, if $|A_{1}|=|A_{2}|$ then we have
\begin{align*}
|A+4\cdot
A|&=|X_{1}+A|+|X_{2}+A|\\&\geq|X_{1}+A_{1}|+|\triangle_{11}|+|X_{2}+A_{2}|+|\triangle_{22}|\\
&\geq5|A|-4-4+2\geq 5|A|-6
\end{align*}
by Lemma \ref{ii}. In the case $|A_{1}|>|A_{2}|$, with the help of Lemma \ref{a+kb} we
obtain
\begin{align*}
|A+4\cdot
A|&=|X_{1}+A|+|X_{2}+A|\geq5|A_{1}|-4+|A_{2}|+4|A_{1}|-4\\
&=5|A|+4(|A_{1}|-|A_{2}|)-8\geq 5|A|-6.
\end{align*}

When $|\ \{i:|\hat{X_{i}}|\leq 3\}\ |=1$, it is easy to see that
$$|A+4\cdot A|\geq5(|A_{1}|+|A_{2}|)+|A_{2}|-4-6\geq5|A|-6.$$

In the case
$|\ \{i:|\hat{X_{i}}|\leq 3\}\ |=2$, we have
$$|A+4\cdot
A|\geq5|A_{1}|-6+|A_{2}|+5|A_{2}|-6+|A_{1}|\geq5|A|-6.$$

Combining the above, we have completed the proof of Theorem \ref{k=4}.
\end{proof}

\subsection*{Acknowledgement}
The authors are grateful to the referee for helpful comments.


\end{document}